\documentclass{elsart3-1}

 \usepackage{graphicx}
 \usepackage{epsfig}

\usepackage{amssymb,MnSymbol,cite}

\usepackage[english,francais]{babel}

\newtheorem{theorem}{Theorem}[section]

\newtheorem{e-proposition}[theorem]{Proposition}
\newtheorem{corollary}[theorem]{Corollary}
\newtheorem{e-definition}[theorem]{Definition\rm}

\renewcommand{\theequation}{\arabic{equation}}
\def\og{\leavevmode\raise.3ex\hbox{$\scriptscriptstyle\langle\!\langle$~}}
\def\fg{\leavevmode\raise.3ex\hbox{~$\!\scriptscriptstyle\,\rangle\!\rangle$}}

\def\Aut{{\rm Aut}}

\def\End{{\rm End}}

\def\GL{{\rm GL}}

\def\id{{\rm id}}

\def\Spec{{\rm Spec\,}}

\def\Tr{{\rm Tr}}

\def\A{{\mathbb A}}
\def\B{{\mathbb B}}

\def\F{{\mathbb F}}

\def\N{{\mathbb N}}

\def\Q{{\mathbb Q}}
\def\R{{\mathbb R}}
\def\Z{{\mathbb Z}}

\def\urep{\vartheta}

\def\Tr{{\rm Tr}}

\def\cC{{\mathcal C}}

\def\cF{{\mathcal F}}

\def\cO{{\mathcal O}}

\def\cR{{\mathcal R}}

\def\cT{{\mathcal T}}

\def\qqq{\,,\,~\forall}

\def\sub{{\rm Sub}_\geq}
\def\conv{{\rm Conv}_\geq}

\newcommand{\ie}{{\it i.e.\/}\ }

\newcommand{\cf}{{\it cf.}}

\def\id{{\mbox{Id}}}

\def\End{{\mbox{End}}}

\def\mc{multiplicatively cancellative }

\def\rmax{\R_+^{\rm max}}

\def\fr{{\rm Fr}}

\def\coo{{[0,1]_{\rm max}}}

\def\arith{{(\wnt,\bar \N)}}

\def\tt{\frak{t}}

\def\nt{\N^{\times}}
\def\wnt{{\widehat{\N^{\times}}}}

\def\nbo{{\overline\N\otimes_\B \overline\N}}

\def\wntb{{\widehat{\N^{\times 2}}}}
\def\arithb{{(\wntb,\nbo)}}
\def\arithc{{(\wntb,\conv(\N\times \N))}}
\def\beps{{\bar \N_\epsilon}}
\def\germ{{{\rm Germ}_{\epsilon=0}(\rmax)}}
\def\vsp{\vspace{.05in}}

\newcommand\blfootnote[1]{%
  \begingroup
  \renewcommand\thefootnote{}\footnote{#1}%
  \addtocounter{footnote}{-1}%
  \endgroup
}

\journal{the Acad\'emie des Sciences}
\begin{document}
% place in the next line the header (rubrique) chosen for your article,
% if you know it (you can also have 2, format : Header1/Header2
\centerline{}
\begin{frontmatter}


% \author{Name\thanksref{label2}}
\selectlanguage{english}
\title{The Arithmetic Site\\\vspace{.05in}
Le Site Arithm\' etique}

% The [label1] can be suppressed if there is only one address for all authors

\selectlanguage{english}
\author[authorlabel1]{Alain Connes},
%\thanks[label1]{2.~wish to thank Ohio State University in which the paper was written.}
\ead{alain@connes.org}
\author[authorlabel2]{Caterina Consani\thanksref{label2}}
\thanks[label2]{Partially supported by the NSF grant DMS 1069218.
}
\ead{kc@math.jhu.edu}

\address[authorlabel1]{Coll\`ege de France,
3 rue d'Ulm, Paris F-75005 France;
I.H.E.S. and Ohio State University}
\address[authorlabel2]{The Johns Hopkins
University Baltimore, MD 21218 USA}

% If you know the dates of reception, and acceptation you can put them now;
%  idem the name of the person presenting the Note

%\medskip
%\begin{center}
%{\small Received *****; accepted after revision +++++\\
%Presented by £££££}
%\end{center}

\begin{abstract}
\selectlanguage{english}

We show that the non-commutative geometric approach to the Riemann zeta function has an algebraic geometric incarnation: the ``Arithmetic Site". This site involves the tropical semiring $\bar\N$ viewed as a sheaf on the topos $\wnt$ dual to the multiplicative semigroup of positive integers. We realize the Frobenius correspondences in the square of the ``Arithmetic Site". 

%{\it To cite this article: A.
%Name1, A. Name2, C. R. Acad. Sci. Paris, Ser. I 340 (2005).}

\vskip 0.5\baselineskip

\selectlanguage{francais}

\noindent{\bf R\'esum\'e} \vskip 0.5\baselineskip \noindent

Le ``Site Arithm\' etique" est  l'incarnation en g\'eom\'etrie alg\'ebrique de l'espace non-commutatif, de nature ad\' elique, qui permet d'obtenir la fonction z\^{e}ta de Riemann comme fonction de d\' enombrement de Hasse-Weil. 
Ce site est construit \` a partir du semi-anneau tropical $\bar\N$ vu comme un faisceau sur le topos $\wnt$ dual du semigroupe multiplicatif des entiers positifs. Nous r\'ealisons les correspondances de Frobenius dans le carr\' e du ``Site Arithm\' etique".

%{\it Pour citer cet article~: A. Name1, A. Name2, C. R. Acad. Sci.
%Paris, Ser. I 340 (2005).}

\end{abstract}
\end{frontmatter}

% now the Version française abrégée, if it exists
%\selectlanguage{francais}
%\section*{Version fran\c{c}aise abr\'eg\'ee}
% Text of your Version française abrégée here.
% Note you do not need to repeat here equations that you use in the
% main text - for example 'voir (3)' is quite acceptable.

\selectlanguage{english}
% main text
\section{Introduction}
\label{}
\blfootnote{Both authors are grateful to Ohio State University where this paper was written.}
\blfootnote{{\it Keywords:}~Site Arithm\'etique, classes d'Ad\`eles, topos, correspondances de Frobenius, caract\'eristique 1.
Arithmetic Site, Ad\`ele class space, topos, Frobenius correspondences, characteristic 1.}

We unveil  the ``Arithmetic Site"  as a ringed topos deeply related to the non-commutative geometric approach to RH. The topos  is the presheaf topos $\wnt$ of functors from the multiplicative semigroup $\N^\times$ of positive integers to the category  of sets. The   structure sheaf is a sheaf of semirings of characteristic $1$ and  (as an object of the topos) is the tropical semiring  $\bar \N:=(\N\cup \infty, \inf,+)$, $\N=\Z_{\geq 0}$,  on which the semigroup $\N^\times$ acts by multiplication.
We prove that  the set of points of the arithmetic site $\arith$ over the maximal compact subring $\coo\subset \rmax$ of the tropical semifield   is the non-commutative space $\Q^\times\backslash \A_\Q/ \hat\Z^*$ quotient of the ad\` ele class space of $\Q$ by the action of the maximal compact subgroup $\hat\Z^*$ of the idele class group. 
In \cite{CC1,CC2}  it was shown that the action of $\R_+^*$ on $\Q^\times\backslash \A_\Q/ \hat\Z^*$  yields the counting distribution whose Hasse-Weil zeta function  is the complete Riemann zeta function. This result is now applied to the arithmetic site to show that its Hasse-Weil zeta function is the complete Riemann zeta function.  The action of $\R_+^*$ on $\Q^\times\backslash \A_\Q/ \hat\Z^*$ indeed corresponds to the action of the Frobenius automorphisms  $\fr_\lambda\in \Aut(\rmax)$, $\lambda\in \R_+^*$, on the  points of $\arith$ over $\coo\subset\rmax$. The square of the arithmetic site over the semifield $\B=(\{0,1\},\max,\times)$ has an unreduced and reduced version. In both cases the underlying topos is $\wntb$. The structure sheaf  in the unreduced case is $\bar\N\otimes_\B \bar\N$ and in the reduced case is the  \mc semiring canonically associated to $\bar\N\otimes_\B \bar\N$.  We determine this latter semiring as the semiring $\conv(\N\times \N)$ of Newton polygons with the operations of convex hull of the union and sum. On both versions there is a canonical action of $\N^{\times 2}$ by endomorphisms $\fr_{n,m}$. By composing this action with the diagonal (given by the product $\mu$) one obtains the Frobenius correspondences $\Psi(\lambda)=\mu\circ \fr_{n,m}$ for rational values  $\lambda=n/m$. 
The Frobenius correspondences $\Psi(\lambda)$ for arbitrary positive real numbers $\lambda$  are 
 realized as curves in the square obtained from the rational case using diophantine approximation. Finally we determine the composition law of these correspondences and show that it is given by the product law in $\R_+^*$ with a subtle nuance in the case of  two irrational numbers whose product is rational. 

This note provides the algebraic geometric space underlying the non-commutative approach to RH. It gives a geometric framework  reasonably  suitable to transpose the conceptual understanding  of the Weil proof in finite characteristic  as in \cite{Gr}.  This translation would require in particular an adequate version of the Riemann-Roch theorem in characteristic $1$.

\section{The arithmetic site}
Given a small category $\cC$ we denote by $\hat \cC$ the topos of contravariant functors from $\cC$ to the category of sets. We let $\nt$ be the category with a single object $*$, $\End(*)=\nt$.
\vsp
\begin{e-definition}\label{site} We define the {\em arithmetic site} $\arith$ as the topos $\wnt$
endowed with the  {\em structure sheaf} $\bar \N:=(\N\cup \infty, \inf,+)$ viewed as a semiring in the topos.
\end{e-definition}
\vsp
Notice that $\wnt\simeq \text{Sh}(\nt,\text{J})$, where $\text{J}$ is the chaotic topology on $\nt$ (\!\!\cite{AGV} Expos\'e IV, 2.6).

\subsection{The points of the topos $\wnt$} 

A point of a topos $\cT$ is defined as a geometric morphism from the topos of sets to $\cT$ (\!\!\cite{AGV,MM}).\vsp

\begin{theorem}\label{thmpointswnt}
$(i)$~The category of points of the topos $\widehat{ \N^\times} $ is canonically equivalent to the category of totally ordered groups isomorphic to non-trivial subgroups of   $(\Q,\Q_+)$, and injective morphisms of ordered groups.

$(ii)$~Let  $\A_f$ be the ring of finite ad\`eles of $\Q$. The space of isomorphism classes of points of $\wnt$ is canonically isomorphic to the double quotient
$
\Q_+^\times\backslash\A_f/\hat\Z^*
$
where $\Q_+^\times$ acts by multiplication on $\A_f$.

\end{theorem}
\vsp
We denote by $\F=\Z_{\rm max}$ the semifield of fractions of the semiring $\bar\N$.\vsp

\begin{corollary}\label{pointswnt} The category of points of the topos $\widehat{ \N^\times} $ is equivalent to the category  of algebraic extensions of the semi-field  $\F=\Z_{\rm max}$ \ie of extensions: 
$
\F\subset K\subset \bar \F=\Q_{\rm max}.
$
The morphisms are the {\em injective} morphisms of semifields.
\end{corollary}

\subsection{The structure sheaf $\bar\N$} 

The next result provides an explicit description of the semiring structure inherited  automatically by  the stalks   
of the sheaf $\bar \N$ on the topos $\wnt$.\vsp

\begin{theorem} \label{structure2} At the point of the topos $\wnt$ associated to the intermediate semifield $\F\subset K\subset \bar \F=\Q_{\rm max}$ the stalk of the structure sheaf $\cO:=\bar \N$  is the semiring $\cO_K:=\{r\in K\mid r\vee 1=1\}$ where $\vee$ denotes addition.
\end{theorem}

\subsection{The points of the arithmetic site $\arith$ over $\coo$}

The following definition provides the notion of point of the arithmetic site over a local semiring.\vsp

\begin{e-definition} Let $R$ be a local semiring. Then a morphism $f:\Spec(R)\to\arith $ is a pair of a point $p$ of $\wnt$ and a local morphism  of semirings $f^\#_p:\cO_p\to R$.
\end{e-definition}
\vsp
The next crucial  statement determines the interpretation of the  space underlying the non-commutative geometric approach to RH in terms of algebraic geometry.\vsp

\begin{theorem} \label{structure3} The points of the arithmetic site  $\arith$ over the maximal compact subring $\coo\subset \rmax$ of the tropical semifield form the quotient $\Q^\times\backslash \A_\Q/ \hat\Z^*$ of the ad\`ele class space of $\Q$ by the action of $\hat\Z^*$. The action of the Frobenius
automorphisms $\fr_\lambda\in \Aut(\coo)$ on these points corresponds to the action of the id\`ele class group (mod $ \hat\Z^*$) on the above quotient of the ad\`ele class space.
\end{theorem}
\vsp
Notice that the quotient $\Q^\times\backslash \A_\Q/ \hat\Z^*$ is the disjoint union of the following two spaces:

$(i)$~$
\Q_+^\times\backslash\A_f/\hat\Z^*
$ is the space of  ad\`ele classes whose archimedean component vanishes. The corresponding points of the arithmetic site $\arith$ are those  which are defined over $\B$;  they are given by the  points of $\wnt$ (Theorem \ref{thmpointswnt}).

$(ii)$~$\Q_+^\times\backslash\left(( \A_f/\hat\Z^*)\times \R_+^*\right)$  is the space of  ad\`ele classes whose archimedean component does not vanish. It is in canonical bijection with rank one subgroups  of $\R$ through the map 
$$
(a,\lambda) \mapsto  \lambda H_a\qqq a\in \A_f/\hat\Z^*, \, \lambda \in  \R_+^*, \  \  H_a:=\{q\in\Q\mid qa\in \hat\Z \}.
$$
%\begin{table}\label{table:1}
%  \epsfig{file=tomato2.ps}
%  \caption{Differential responding of tomatic stimulation in the
 %   brain at different frequencies.}
%\end{table}

\subsection{Hasse-Weil formula for the Riemann zeta function}

In order to count the number of fixed points of the Frobenius action on points of $\arith$ over $\coo$ we let $
\urep_u\xi(x)=\xi(u^{-1}x)
$ be the scaling action of the id\`ele class group $G=\GL_1(\A_\Q)/\GL_1(\Q)$  on functions on the ad\`ele class space $\A_\Q/\Q^*$ and use the trace formula (\!\!\cite{Co-zeta,Meyer,CMbook})
 in the form ($\Sigma_\Q=$ places of $\Q$, $d^*u$ multiplicative Haar measure)
\begin{equation}\label{trace}
\Tr_{\rm distr}\left(\int_G h(u)\urep_ud^*u\right )=\sum_{v\in\Sigma_\Q}\int_{\Q^*_v}\,\frac{h(u^{-1})}{|1-u|}\,d^*u.
\end{equation}
We apply \eqref{trace}  to test functions of the form $h(u)=g(|u|)$ where the support of  $g$ is contained in
$(1,\infty)$ and $\vert u\vert$ is the module. The invariance of $h$ under the kernel $\hat\Z^*$ of the module $G\to \R_+^*$ corresponds at the geometric level  to taking the quotient of the ad\`ele class space by the action of $\hat\Z^*$. Using  Theorem \ref{structure3} and \cite{CC2}, \S 2, one obtains  the  counting distribution $N(u)$, $u\in [1,\infty)$ associated to the Frobenius action on points of $\arith$ over $\coo$.\vsp

\begin{theorem}\label{mainthm} The zeta function $\zeta_N$ associated by the equation
\begin{equation}\label{logzetabis}
    \frac{\partial_s\zeta_N(s)}{\zeta_N(s)}=-\int_1^\infty  N(u)\, u^{-s}d^*u
\end{equation}
to the counting distribution $N(u)$ is the  {\em complete} Riemann zeta function  $\zeta_\Q(s)=\pi^{-s/2}\Gamma(s/2)\zeta(s)$.
\end{theorem}
\vsp
In \cite{CC1}, equation \eqref{logzetabis} was shown (following a suggestion made in \cite{Soule})  to be  the limit, when $q\to 1$, of the Hasse-Weil formula for counting functions over finite fields $\F_q$. 

\section{The square of the arithmetic site}

\subsection{The  unreduced square $\arithb$}

Given a partially ordered set $J$,  we let $\sub(J)$ be the set of subsets $E\subset J$ which are hereditary, \ie such that $x\in E\implies y\in E$, $\forall y\geq x$.  Then  $\sub(J)$ endowed with the operation $E \oplus E':=E\cup E'$ is a $\B$-module. We refer to \cite{PR} for the general treatment of tensor products of semi-modules.\vsp

\begin{e-proposition}\label{tensring} 
$(i)$~Let $\N\times \N$ be endowed with the partial order $ (a,b)\leq (c,d)\iff a\leq c \  \& \  b\leq d$. Then one has a canonical isomorphism of $\B$-modules $\overline\N\otimes_\B \overline\N\simeq \sub(\N\times \N).
$

\noindent $(ii)$~There exists on the $\B$-module  $S=\overline\N\otimes_\B \overline\N$ a unique bilinear multiplication
such that, using multiplicative notation where $q$ is a formal variable,  one has
\begin{equation}\label{square3}
(q^ a\otimes_\B q^b)(q^ c\otimes_\B q^d)=q^ {a+c}\otimes_\B q^{b+d}.
\end{equation} 
$(iii)$~The  multiplication \eqref{square3} turns $\overline\N\otimes_\B \overline\N$ into a semiring of characteristic $1$.\newline
$(iv)$~The following formula defines an action of $\N^\times\times \N^\times$  by endomorphisms on $\overline\N\otimes_\B \overline\N$
\begin{equation*}\label{bifrob}
\fr_{n,m}(\sum q^a\otimes_\B q^b):=\sum q^{na}\otimes_\B q^{mb}.
\end{equation*}
\end{e-proposition}\vsp

\begin{e-definition}\label{site2} The {\em unreduced square} $\arithb$ of the arithmetic site $\arith$ is the topos $\wntb$ with the  {\em structure sheaf} $\nbo$, viewed as a semiring in the topos.
\end{e-definition}

\subsection{The Frobenius correspondences}\label{sectfrob}
The product in the semiring $\overline \N$ yields a morphism of semirings $\mu:
(\overline\N\otimes_\B \overline\N)\to \overline \N$, given on simple tensors by 
$
\mu(q^ a\otimes_\B q^b)=q^{a+b}.
$
\vsp

\begin{e-proposition}\label{frobcor} 
$(i)$~The range of the morphism $\mu\circ \fr_{n,m}: \overline\N\otimes_\B \overline\N\to \overline\N$ only depends, up to canonical isomorphism, on the ratio $r=n/m$. Assuming that $n,m$ are relatively prime, this range contains the ideal
$$
\{q^a\mid a\geq (n-1)(m-1)\}\subset \bar\N.
$$
$(ii)$~Let $r=n/m$, $q\in(0,1)$  and let $m_r: \overline\N\otimes_\B \overline\N\to \rmax$ be given by
$$
  m_r\left(\sum (q^{n_i}\otimes_\B q^{m_i})\right)=q^\alpha, \   \alpha= \inf (r n_i+ m_i).
$$
Up to canonical isomorphism of their ranges, the  morphisms $\mu\circ \fr_{n,m}$ and $m_r$ are equal.
\end{e-proposition}
\vsp
Proposition \ref{frobcor} $(ii)$~allows one to extend the definition of the Frobenius correspondence  to  arbitrary positive real numbers. \vsp

\begin{e-proposition}\label{frobcorres} 
$(i)$~Let $\lambda\in \R_+^*$ and $q\in (0,1)$ then the following formula defines a homomorphism\begin{equation}\label{frobc1}
\cF(\lambda,q):\overline\N\otimes_\B \overline\N\to\rmax, \ \  \cF(\lambda,q)\left(\sum (q^{n_i}\otimes_\B q^{m_i})\right)=q^\alpha, \   \alpha= \inf (\lambda n_i+m_i).
\end{equation} 
$(ii)$~The  semiring   $\cR(\lambda):={\rm Im}( \cF(\lambda,q))$ is independent, up to canonical isomorphism, of  $q\in (0,1)$.\newline
$(iii)$~The semirings $\cR(\lambda)$ and $\cR(\lambda')$ are isomorphic if and only if $\lambda'=\lambda$ or $\lambda'=1/\lambda$.
\end{e-proposition}

\subsection{The reduced square $\arithc$}

Let $R$ be a semiring without zero divisors and $\iota:R\to {\rm Frac}R$ the canonical morphism to the semifield of fractions. It is not true in general that $\iota$ is injective  (\cf \cite{DC}). We shall refer to $\iota(R)$  as the {\em reduced} semiring of $R$.\vsp

\begin{e-definition} We let $\conv(\N\times \N)$ be the set of closed convex subsets $C$ of the quadrant $Q:=\R_+\times \R_+$ such that $(i)$~$C+Q=C$ and $(ii)$~the extreme points $\partial C$ belong to $\N\times \N\subset Q$.
\end{e-definition}
\vsp
The set $\conv(\N\times \N)$ is a semiring for the operations of convex hull of the union and sum.\vsp

\begin{e-proposition}\label{conv}
$(i)$~The semiring $\conv(\N\times \N)$ is multiplicatively cancellative.\newline
$(ii)$~The homomorphism $\gamma:\nbo\simeq\sub(\N\times \N)\to\conv(\N\times \N)$ given by convex hull is the same as the homomorphism  $\iota:\nbo\to \iota(\nbo)$. \newline
$(iii)$~
Let $R$ be a   multiplicatively cancellative semiring and $\rho:\nbo\to R$
a homomorphism of semirings such that $\rho^{-1}(\{0\})=\{0\}$. Then there exists a unique semiring homomorphism $\rho':\conv(\N\times \N)\to R$ such that 
$\rho=\rho'\circ \gamma$.
\end{e-proposition}\vsp

\begin{e-definition}\label{site3} The {\em reduced square} $\arithc$ of the arithmetic site $\arith$ is the topos $\wntb$ with the  {\em structure sheaf} $\conv(\N\times \N)$, viewed as a semiring in the topos.
\end{e-definition}

\begin{figure}
\begin{center}
\includegraphics[scale=0.6]{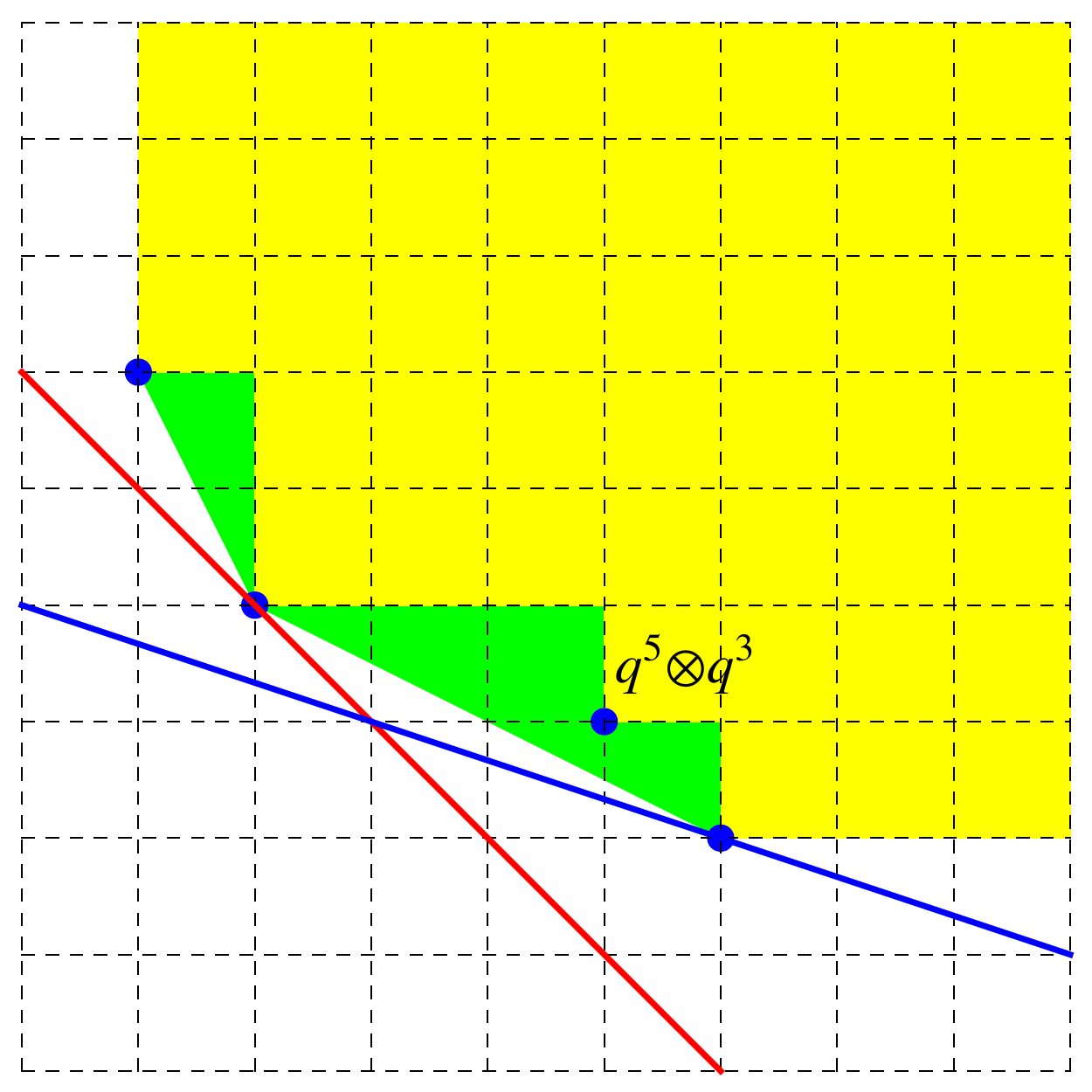}
\caption{Typical element $E\in \sub(\N\times \N)$ (in yellow); its image under $\mu$ (red), under the Frobenius correspondence for $\lambda=\frac 13$ (blue) and under $\gamma$ (green or yellow). Removing the element $q^5\otimes_\B q^3$ does not alter  the convex hull $\gamma(E)$.}\label{adsub1}
\end{center}
\end{figure}
\section{Composition of Frobenius correspondences}\label{sectcomp}

\subsection{Reduced correspondences}

\begin{e-definition}\label{defcorr} A {\em reduced correspondence} over the arithmetic site $\arith$ is given by a triple $(R,\ell,r)$ where $R$ is a \mc semiring, $\ell, r: \bar\N\to R$ are semiring morphisms such that $\ell^{-1}(\{0\})=\{0\}$, $r^{-1}(\{0\})=\{0\}$
and that $R$ is generated by $\ell(\bar\N)r(\bar\N)$. 
\end{e-definition}
By construction, \cf~Proposition \ref{frobcorres},  the Frobenius correspondence 
gives a reduced correspondence: 
\begin{equation}\label{frobcordef}
\Psi(\lambda):=(R,\ell(\lambda),r(\lambda)), \  \  R:=\cR(\lambda), \  \  \ell(\lambda)(q^n):=\cF(\lambda,q)(q^n\otimes 1), \  \  r(\lambda):=\cF(\lambda,q)(1\otimes q^n).
\end{equation}
By \eqref{frobc1} one gets that the elements of $\cR(\lambda)$ are powers $q^\alpha$ where $\alpha \in \N+\lambda \N$ and that the morphisms $\ell(\lambda)$ and $r(\lambda)$ are described as follows:
\begin{equation*}\label{morphism3}
\ell(\lambda)(q^n)q^\alpha=q^{\alpha +n\lambda}, \    \  r(\lambda)(q^n)q^\alpha=q^{\alpha +n}
\end{equation*}

\subsection{The composition of the correspondences $\Psi(\lambda)\circ \Psi(\lambda')$}

The composition $\Psi(\lambda)\circ \Psi(\lambda')$ of the Frobenius correspondences is obtained as the left and right action of $\bar \N$ on the reduced semiring of the tensor product $\cR(\lambda)\otimes_{\bar\N} \cR(\lambda')$. In order to state the general result we introduce a variant  $\id_\epsilon$ of the identity correspondence.
We let $\germ$ be the semiring of  germs of continuous functions from a neighborhood of $0\in \R$ to $\rmax$, endowed with the pointwise operations. 
Let $\beps$ be the the sub-semiring of  $\germ$ generated, for fixed $q\in (0,1)$, by $q$ and $\fr_{1+\epsilon}(q)=q^{1+\epsilon}$. $\beps$ is independent, up to canonical isomorphism, of the choice of $q\in (0,1)$.\vsp

\begin{e-definition}\label{defndef} The tangential deformation of the identity correspondence is given by the triple $(\beps,\ell_\epsilon,r_\epsilon)$ where  $\ell_\epsilon(q^n):=\fr_{1+\epsilon}(q^n)$ and $r_\epsilon(q^n):=q^n$, $\forall n\in \N$.
\end{e-definition}\vsp

\begin{theorem}\label{thmcomp} Let $\lambda, \lambda' \in \R_+^*$ such that $\lambda\lambda'\notin \Q$. The composition of the Frobenius correspondences is then given by 
$$
\Psi(\lambda)\circ \Psi(\lambda')=\Psi(\lambda\lambda')
$$
The same equality holds if $\lambda$ and $\lambda'$ are rational. When  $\lambda, \lambda'$ are irrational and $\lambda\lambda'\in \Q$,
$$
\Psi(\lambda)\circ \Psi(\lambda')=\Psi(\lambda\lambda')\circ \id_\epsilon=\id_\epsilon\circ \Psi(\lambda\lambda')
$$
where $\id_\epsilon$ is the tangential deformation of the identity correspondence.
\end{theorem}

%\begin{table}\label{table:1}
%\epsfig{sqcurve8.eps}  \caption{Differential responding of tomatic stimulation in the
%   brain at different frequencies.}
%\end{table}

%\section{Thanks}
%{\it } \cite{b15}

% etc, etc

% \appendix

% \section{}
% \label{}

% The Acknowledgements are an un-numbered section
%\section*{Acknowledgements}
% Acknowledgements text here

\end{document}